\documentclass[11pt,leqno]{amsart}
\usepackage{amssymb,booktabs,array,microtype,needspace}
\usepackage[hidelinks]{hyperref}

\newtheorem{theorem}{Theorem}[section]
\newtheorem{lemma}[theorem]{Lemma}
\newtheorem{proposition}[theorem]{Proposition}
\theoremstyle{definition}
\newtheorem{definition}[theorem]{Definition}
\numberwithin{equation}{section}
\DeclareMathOperator{\height}{height}
\DeclareMathOperator{\Hom}{hom}
\DeclareMathOperator{\Inj}{inj}
\DeclareMathOperator{\comp}{comp}
\DeclareMathOperator{\Av}{Av}
\DeclareMathOperator{\Top}{Top}
\DeclareMathOperator{\MaxInt}{MaxInt}
\newcommand{\one}{\mathbf1}
\newcommand{\Q}{\mathbb Q}
\newcommand{\R}{\mathcal R}
\newcommand{\Rp}{\mathcal R_+}
\newcommand{\dual}{^{\partial}}
\newcommand{\M}{\mathrm M}
\newcommand{\J}{\mathrm J}
\newcommand{\Js}{\mathrm J^{*}}
\newcommand{\Log}{\mathrm L}
\newcommand{\setm}[2]{\{#1\mid #2\}}

\hypersetup{pdftitle={Two Short Proofs of the A-to-A Theorem},
 pdfauthor={George Gr\"atzer}}

\begin{document}
\title{Two Short Proofs of the $A^A$ Theorem}
\author{George Gr\"atzer}
\address{Department of Mathematics, University of Manitoba,
  Winnipeg, MB R3T 2N2, Canada}
\email{gratzer@mac.com}
\date{}
\subjclass[2020]{06A06}
\keywords{Finite ordered set, order-preserving map, exponentiation,
  reconstruction, homomorphism count}

\begin{abstract}
In 1984, D.~Duffus asked whether $A^A\cong B^B$ implies $A\cong B$
for finite, nonempty ordered sets $A$ and~$B$. J.\,D. Farley and the
author proved that it does. 

In this paper, we give two short proofs using all intervals
of maximum height and the elements with exactly one lower cover.
The first proof averages over these intervals; the second takes their
disjoint sum and constructs three invariants. Both avoid the auxiliary
ordering, direct factorization, and canonical coordinates of the earlier
proof. 
\end{abstract}
\maketitle

\section{Introduction}\label{S:intro}

\subsection{The background}
The $A^A$ problem is easy to state. For finite, nonempty ordered sets
$X$ and~$Y$, let $X^Y$ denote the set of order-preserving
maps $Y\to X$, ordered pointwise. In 1984, D.~Duffus
\cite[p.~90]{DufHD} asked whether $A^A$ determines~$A$.

Recently, J.\,D. Farley and I answered the question affirmatively, see \cite{FG}. 

\begin{theorem}\label{T:main}
Let $A$ and $B$ be finite, nonempty ordered sets. If $A^A\cong B^B$,
then $A\cong B$.
\end{theorem}

\subsection{Shorter proofs}
The proof was fairly long. Somewhat later, I succeeded in finding  a~shorter proof, see Section~\ref{S:average}.

I realized that I had a~unique opportunity
to test how competitive ChatGPT was with my own mathematical
abilities. I decide to challenge ChatGPT to find another short proof.

Section~\ref{S:invariants} presents the response of ChatGPT. 

\subsection{Outline}
The preliminaries shared by the two short
proofs have been collected in Section~\ref{S:prelim}.
My short proof is in Section~\ref{S:average}.
ChatGPT's proof is in Section \ref{S:invariants}.

Finally, Section \ref{S:compare} discusses how different techniques are used in the three proofs. 

\subsection{What this challenge means?} My conclusion from this experiment
is that ChatGPT is fully competitive with my own mathematical
research abilities. If I look just a~year ahead and consider how
much further ChatGPT may develop, I expect the comparison will be
overwhelmingly in ChatGPT's favor.

\section{Common preliminaries}\label{S:prelim}

Unless stated otherwise, ordered sets are finite and nonempty.
Auxiliary disjoint sums and the constructions $\J$ and $\Js$ below may
be empty. The empty ordered set is the zero for disjoint sums and is
absorbing for direct products. Heights and intervals of maximum height
are used only for nonempty ordered sets.
An ordered set is \emph{connected} if its comparability graph is
connected, \emph{bounded} if it has a~least and a~greatest element,
and \emph{nontrivial} if it has more than one element.
We write $\one$ for the one-element ordered set.
The nonempty hypothesis in Theorem~\ref{T:main} is necessary:
$\varnothing^{\varnothing}\cong\one^{\one}\cong\one$.

\subsection{Height and intervals}
We write $X+Y$ for the disjoint sum, with no comparabilities between
the summands, $nX$ for the sum of $n$ copies of~$X$, $X\times Y$ or $XY$
for the direct product, and $X\dual$ for the dual. We allow $n=0$, with
$0X=\varnothing$. We write $\comp(X)$ for the number of connected
components of $X$, and $\height(X)$ for the maximum of $|C|-1$ over
all chains~$C$ of~$X$.

For $a\leq b$, the interval $[a,b]=\setm{x}{a\leq x\leq b}$ is
bounded and connected. Put
\[
 \MaxInt(X)=\setm{[a,b]\subseteq X}{\height([a,b])=\height(X)}.
\]
This set is nonempty: the endpoints of a~longest chain give such an
interval. We count distinct intervals separately, even when they are
isomorphic. An isomorphism of ordered sets induces a~bijection between
their maximum-height intervals.

\begin{lemma}\label{L:height}
For ordered sets $X$ and $Y$,
\begin{align}\label{E:heights}
 \height(X\times Y)&=\height(X)+\height(Y),\\
 \height(X^Y)&=|Y|\height(X).
\end{align}
\end{lemma}
\begin{proof}
Every strict step in a~chain in $X\times Y$ increases at least one
coordinate. The total number of such increases is at most
$\height(X)+\height(Y)$. Following longest chains in the two
coordinates one after the other attains this bound.

In a~chain of maps $Y\to X$, each value can increase strictly at most
$h=\height(X)$ times. Every strict step changes at least one value,
so $\height(X^Y)\leq |Y|h$. For the reverse inequality, choose a~chain
$x_0<\cdots<x_h$ in $X$ and a~linear extension of~$Y$.
Start with the constant map $x_0$. For each $j=1,\dots,h$, change
the values from $x_{j-1}$ to $x_j$ one at a~time, in reverse order of
the linear extension. The elements already changed form an up-set,
so every intermediate map is order-preserving. There are $|Y|h$
strict steps.
\end{proof}

The following description of intervals is also used in~\cite{FG}.
We state the correspondence as a~bijection because its multiplicities
matter in both proofs.

\begin{lemma}\label{L:intervals}
Let $Y_1,\dots,Y_r$ be the connected components of~$Y$, and choose
$y_i\in Y_i$ for each~$i$. For every $[f,g]\in\MaxInt(X^Y)$, the maps
$f$ and $g$ are constant on each~$Y_i$, and
\[
 [f,g]\cong\prod_{i=1}^r[f(y_i),g(y_i)]^{Y_i}.
\]
The rule
$[f,g]\mapsto\bigl([f(y_1),g(y_1)],\dots,[f(y_r),g(y_r)]\bigr)$
is a~bijection from $\MaxInt(X^Y)$ onto $\MaxInt(X)^r$.
\end{lemma}
\begin{proof}
Let $h=\height(X)$ and $[f,g]\in\MaxInt(X^Y)$. Along a~chain in
$[f,g]$, the value at $y$ can increase at most
$\height([f(y),g(y)])$ times. Since every strict step changes a~value,
\[
 |Y|h=\height([f,g])
 \leq\sum_{y\in Y}\height([f(y),g(y)])\leq |Y|h.
\]
Thus every interval $[f(y),g(y)]$ has height~$h$. A~longest chain in
this interval contains both endpoints and is a~longest chain in~$X$.
An element below $f(y)$ or above $g(y)$ would extend it. Hence $f(y)$
is minimal and $g(y)$ is maximal in~$X$. Comparable minimal elements
are equal, as are comparable maximal elements. Isotonicity and
connectedness therefore make $f$ and $g$ constant on each component.
Restriction to the components gives the displayed order isomorphism.

Conversely, choose an interval of maximum height for each component.
The maps constant on the components at the chosen endpoints span
an interval of height $\sum_i|Y_i|h=|Y|h$. They are uniquely determined
by the chosen intervals, and their endpoints recover those intervals.
This proves the bijection.
\end{proof}

Every interval in $XY$ has the form $I\times J$, and the height formula
shows that it has maximum height exactly when $I\in\MaxInt(X)$ and
$J\in\MaxInt(Y)$. Thus these pairs also give a~bijection of actual
intervals, with multiplicities retained.

\Needspace{7\baselineskip}
\begin{lemma}\label{L:start}
Let $A^A\cong B^B$. If $A$ is an antichain, then $A\cong B$.
If $\height(A)\geq1$, then $\height(B)\geq1$ and
$|A|\height(A)=|B|\height(B)$.
\end{lemma}
\begin{proof}
The constant maps embed $X$ in $X^X$. If $A$ is an antichain, then
$A^A$ is an antichain with $|A|^{|A|}$ elements. Therefore $B^B$,
and then $B$, are antichains. We obtain $|A|^{|A|}=|B|^{|B|}$.
Since $n\mapsto n^n$ is strictly increasing on the positive integers,
the antichains have the same size.

If $\height(A)\geq1$, then $A^A\cong B^B$ has positive height,
so $B$ is not an antichain. Equation~\eqref{E:heights} gives the
required equality.
\end{proof}

\subsection{Homomorphism counts}
For ordered sets $Z$ and $X$, let $\Hom(Z,X)$ be the number of
order-preserving maps $Z\to X$. For connected~$Z$,
\begin{align}
 \Hom(Z,X+Y)&=\Hom(Z,X)+\Hom(Z,Y),\label{E:hom1}\\
 \Hom(Z,XY)&=\Hom(Z,X)\Hom(Z,Y),\label{E:hom2}\\
 \Hom(Z,X\dual)&=\Hom(Z\dual,X),\label{E:hom3}\\
 \Hom(\one,X)&=|X|,\label{E:hom4}\\
  \Hom(Z,X)&\geq|X|>0.\label{E:hom5}
\end{align}
For the first identity, the image of a~connected order lies in one
summand. The next two follow by taking coordinate maps and reversing
both orders, respectively. Constant maps give the last inequality.
When the target is empty, the count is zero for nonempty~$Z$;
the algebraic identities still hold, but the strict inequality does not.

We use the ordered-set form of Lov\'asz's homomorphism-count
theorem~\cite{LovFG}. The following proof is included for completeness.

\begin{lemma}\label{L:lovasz}
Let $X$ and $Y$ be finite ordered sets, possibly empty.
If \[\Hom(Z,X)=\Hom(Z,Y)\] for every finite, nonempty connected ordered
set $Z$, then $X\cong Y$.
\end{lemma}
\begin{proof}
A~map from a~disjoint sum is determined by independent restrictions
to its components, so the equality extends to every finite, nonempty
ordered set~$Z$.

Write $\Inj(Z,X)$ for the number of injective order-preserving maps.
For a~partition $\pi$ of $Z$, generate a~relation on its blocks from
the inequalities of $Z$ and take its reflexive transitive closure.
Call $\pi$ admissible if this relation is antisymmetric, and write
$Z/\pi$ for the resulting ordered set. The kernel partition of an
order-preserving map to $X$ is admissible: a~cycle among distinct blocks
would force their images to be equal. Conversely, an injective
order-preserving map $Z/\pi\to X$ induces a~map with kernel exactly~$\pi$. Consequently,
\[
 \Hom(Z,X)=\sum_{\pi\text{ admissible}}\Inj(Z/\pi,X).
\]
The discrete partition contributes $\Inj(Z,X)$; every other admissible
partition gives a~smaller ordered set, which may be disconnected.
We prove $\Inj(Z,X)=\Inj(Z,Y)$ by induction on $|Z|$, for all
nonempty~$Z$ at once. Assuming it for all ordered sets smaller
than~$Z$, the displayed formula and its counterpart for $Y$ agree
term by term except possibly in the discrete term, and their left
sides agree by the first paragraph. Hence $\Inj(Z,X)=\Inj(Z,Y)$.

Taking $Z=\one$ first, gives $|X|=|Y|$. If these are zero, the result
is immediate. Otherwise, the identity maps and the equality of injective
counts give injective order-preserving maps $f\colon X\to Y$ and $g\colon Y\to X$.
Both are bijections. The order-preserving permutation $gf$ of the finite
set $X$ has an order-preserving inverse, since its inverse is a~power
of $gf$. Thus $f^{-1}=(gf)^{-1}g$ is order-preserving, and $f$ is
an order isomorphism.
\end{proof}

One can also obtain Lemma~\ref{L:lovasz} directly from the theorem for
finite relational structures. For a~finite binary relation, take its
reflexive transitive closure and identify elements related in both
directions. The quotient is an ordered set. Every relation-preserving
map into an ordered set factors uniquely through this quotient.
The empty test structure has exactly one map into every target.
Thus tests by ordered sets suffice for targets that are ordered sets.

\subsection{Elements with one lower cover}
An element $c$ is a~\emph{lower cover} of the element $x$, if $c<x$ and no element
lies strictly between them. Define $\J(X)$ to be the set of elements
with exactly one lower cover, and $\Js(X)$ to be the set of elements
with exactly one upper cover, in both cases with the inherited order.
These constructions preserve isomorphisms and commute with disjoint
sums, including empty sums. Also,
\[
 \Js(X)=\J(X\dual)\dual.
\]
For $y\in Y$, write $\uparrow y=\setm{z\in Y}{y\leq z}$.

\begin{lemma}\label{L:log}
Let $B$ be bounded. For $j\in\J(B)$ and $y\in Y$, let $f_{j,y}$
take the value $j$ on $\uparrow y$ and $0_B$ elsewhere. Then
$(j,y)\mapsto f_{j,y}$ is an order isomorphism
\begin{equation}\label{E:J-bounded-power}
 \J(B)\times Y\dual\cong\J(B^Y).
\end{equation}
Dually, $\Js(B^Y)\cong\Js(B)\times Y\dual$.
\end{lemma}
\begin{proof}
Each $f_{j,y}$ is order-preserving because $\uparrow y$ is an up-set
and $0_B\leq j$.

Take $f\in\J(B^Y)$. It is not the zero map, so its support
\[S=\setm{x\in Y}{f(x)\ne0_B}\] is a~nonempty up-set.
For a~minimal $x\in S$, all values below $x$ are zero. Every lower
cover $c$ of $f(x)$ therefore gives an order-preserving map obtained
by changing only $f(x)$ to~$c$. This map is a~lower cover of~$f$:
an intermediate map would have to differ at that coordinate alone,
with a~value strictly between $c$ and $f(x)$.
Every nonzero element of a~finite bounded order has a~lower cover.
Distinct choices of $x$ or $c$ give distinct lower covers of~$f$.
Thus $S$ has exactly one minimal element~$y$, and
$j=f(y)\in\J(B)$. Finiteness gives $S=\uparrow y$.

Suppose that $f(z)\ne j$ for some $z>y$, and choose a~minimal such~$z$.
Then $f(w)\in\{0_B,j\}$ for every $w<z$, and $j<f(z)$.
A~maximal element $c$ of $[j,f(z))$ is a~lower cover of $f(z)$ in~$B$
and satisfies $c\geq j$. Changing only $f(z)$ to $c$ preserves order
and gives a~lower cover of~$f$, distinct from the one obtained by
lowering $f(y)$. This contradiction proves that $f=f_{j,y}$.

Conversely, let $j_-$ be the unique lower cover of~$j$. Every element
strictly below $j$ lies below $j_-$: extend it to a~maximal element
below $j$. If $g<f_{j,y}$, then $g(y)<j$, since $g(y)=j$ would force
$g=f_{j,y}$ by isotonicity and $g\leq f_{j,y}$. Hence $g(y)\leq j_-$.
The map obtained from $f_{j,y}$ by changing only its value at $y$
to $j_-$ is order-preserving and dominates every such~$g$.
It is therefore the unique lower cover of $f_{j,y}$.

Finally,
\[
 f_{j,y}\leq f_{j',y'}
 \quad\text{is equialent to}\quad j\leq j'\text{ and }y\geq y'.
\]
For the forward implication, evaluate at $y$; the value on the right
cannot be zero. The reverse implication follows from containment of
the supports. This proves the order isomorphism. For the dual formula,
use $(B^Y)\dual\cong(B\dual)^{Y\dual}$ to obtain
\[
 \Js(B^Y)\cong\bigl(\J(B\dual)\times Y\bigr)\dual
 \cong\Js(B)\times Y\dual.\qedhere
\]
\end{proof}

Call $X$ \emph{componentwise bounded} if every connected component
is bounded.

\begin{lemma}\label{L:J-product}
For componentwise bounded $X$ and $Y$,
\begin{equation}\label{E:J-prod}
 \J(XY)\cong\comp(Y)\J(X)+\comp(X)\J(Y),
\end{equation}
and likewise for $\Js$.
\end{lemma}
\begin{proof}
A~lower cover in a~direct product changes exactly one coordinate,
and changes it along a~cover in that factor. Changing both coordinates,
or changing one by more than a~cover, allows an intermediate element.
In a~finite bounded order the least element is the only element with
no lower cover. Thus, for bounded $B,C$,
\[
 \J(BC)=\bigl(\J(B)\times\{0_C\}\bigr)
       +\bigl(\{0_B\}\times\J(C)\bigr).
\]
The two displayed suborders are incomparable. Summing this identity
over the pairs of components of $X$ and $Y$ gives~\eqref{E:J-prod}.
The upper-cover formula follows by duality.
\end{proof}

\begin{lemma}\label{L:J-height}
If $B$ is bounded and nontrivial, then $\J(B)$ and $\Js(B)$ are
nonempty and have height less than $\height(B)$.
For $B=\one$, both are empty.
\end{lemma}
\begin{proof}
Finiteness and nontriviality give an element covering $0_B$ and
an element covered by $1_B$. These have exactly one lower cover
and exactly one upper cover, respectively. A~chain in $\J(B)$ can
be extended by adding $0_B$ below it; a~chain in $\Js(B)$ can be
extended by adding $1_B$ above it. The singleton case follows
directly from the definitions.
\end{proof}

\section{First proof: averaging over intervals}\label{S:average}

\subsection{An elementary logarithm}
For bounded $B$, put $\Log(B)=\J(B)\dual$.
Lemmas~\ref{L:log}--\ref{L:J-height} give, for bounded $B,C$ and any~$Y$,
\begin{align}
 \Log(BC)&\cong\Log(B)+\Log(C),\label{E:L-product}\\
 \Log(B^Y)&\cong Y\times\Log(B).\label{E:L-power}
\end{align}
For the second identity, dualize~\eqref{E:J-bounded-power}.
Moreover, $\Log(\one)=0$, and for nontrivial bounded $B$,
\begin{equation}\label{E:L-height}
 \Log(B)\ne0,\qquad \height(\Log(B))<\height(B).
\end{equation}

\subsection{Rational combinations and cancellation}
Let $\R$ be the rational vector space with basis the isomorphism types
of nonempty connected ordered sets. Direct products of connected
ordered sets are connected: comparability paths in the two factors
can be followed one coordinate at a~time. Thus direct products extend
bilinearly to a~multiplication on~$\R$.
Let $\Rp$ consist of the finite sums $R=\sum_Kr_KK$ with
$r_K\in\Q_{\geq0}$. Identify an ordered set with the sum of its
connected components, counted with their multiplicities; identify
the empty set with~$0$. Uniqueness of connected components shows that
two ordered sets have the same representation exactly when they are
isomorphic. The components of a~product are the products of components,
so this identification respects multiplication.

For connected $Z$, extend $\Hom(Z,-)$ to~$\R$ by linearity.
Equation~\eqref{E:hom2} shows that it is multiplicative. It is positive
on every nonzero element of~$\Rp$.

\begin{lemma}\label{L:cancel}
For $R,S,T\in\Rp$ with $T\ne0$,
\begin{equation}\label{E:cancel}
 RT=ST\quad\text{implies that}\quad R=S.
\end{equation}
\end{lemma}
\begin{proof}
For connected $Z$, apply $\Hom(Z,-)$ and cancel the positive factor
$\Hom(Z,T)$. This gives $\Hom(Z,R)=\Hom(Z,S)$.
Choose a~common denominator $d\geq1$ for the coefficients of $R$ and
$S$, and represent $dR$ and $dS$ by finite ordered sets $X$ and $Y$,
possibly empty. Their counts agree for every connected $Z$, so
Lemma~\ref{L:lovasz} gives $X\cong Y$. Hence $R=S$.
\end{proof}

For $R\in\Rp - \{0\}$, let $\height(R)$ be the largest height of a~$K$ with
$r_K>0$, and let $\Top(R)$ be the sum of the terms at that height.

\begin{definition}\label{D:Av}
For nonzero $R=\sum_Kr_KK\in\Rp$ of height~$h$, put
\[
 \Av(R)=
 \frac{\displaystyle\sum_{\height(K)=h}r_K
                   \sum_{I\in\MaxInt(K)}\Log(I)}
      {\displaystyle\sum_{\height(K)=h}r_K\,|\MaxInt(K)|}.
\]
\end{definition}

The denominator is positive, and every interval occurring in the
numerator is bounded, so this defines an element of~$\Rp$.
For a~nonempty ordered set~$X$, it gives
\[
 \Av(X)=\frac1{|\MaxInt(X)|}\sum_{I\in\MaxInt(X)}\Log(I).
\]
Clearing denominators writes every nonzero $R\in\Rp$ as $X/n$ for
a~nonempty ordered set $X$ and an integer $n\geq1$. Directly from
the definition,
\begin{align}\label{E:Av-scale}
 \Av(X/n)&=\Av(X),\\ 
 \Av(R)&=\Av(\Top(R)).
\end{align}

\begin{lemma}\label{L:Av-height}
If $\height(R)=0$, then $\Av(R)=0$. If $\height(R)>0$, then
$\Av(R)\ne0$ and $\height(\Av(R))<\height(R)$.
\end{lemma}
\begin{proof}
At height zero every interval in the average is a~singleton, and its
logarithm is zero. At positive height every interval is bounded and
nontrivial, so~\eqref{E:L-height} applies. There are finitely many
terms, all with nonnegative coefficients, and at least one has a~positive
coefficient. Their sum is nonzero and has smaller height.
\end{proof}

\begin{lemma}\label{L:Av-rules}
For nonzero $R,S\in\Rp$ and ordered sets $X,Y$,
\begin{align*}
\Av(RS)&=\Av(R)+\Av(S),\\
 \Av(X^Y)&=Y\Av(X).
\end{align*}
\end{lemma}

\begin{proof}
By~\eqref{E:Av-scale}, the first identity reduces to ordered sets
$R=X'$, $S=Y'$. The maximum-height intervals in $X'Y'$ are exactly
the products $I\times J$, with $I\in\MaxInt(X')$ and
$J\in\MaxInt(Y')$. Average~\eqref{E:L-product} over all these pairs.

For the second identity, put $m=|\MaxInt(X)|$ and write
$Y=Y_1+\cdots+Y_r$ as a~sum of connected components.
By Lemma~\ref{L:intervals}, $\MaxInt(X^Y)$ is in bijection with
$\MaxInt(X)^r$. Equations~\eqref{E:L-product} and~\eqref{E:L-power}
show that the tuple $(I_1,\dots,I_r)$ contributes
$\sum_iY_i\Log(I_i)$. A~given interval $I$ occurs in position $i$
of exactly $m^{r-1}$ tuples. Consequently,
\[
 \Av(X^Y)=\frac1{m^r}\sum_{(I_i)}\sum_iY_i\Log(I_i)
 =\sum_iY_i\cdot\frac1m\sum_I\Log(I)=Y\Av(X).\qedhere
\]
\end{proof}

\begin{lemma}\label{L:Av-inject}
If $R\in\Rp - \{0\}$ and $T=R+\Av(R)$, then
\[
 R=T-\Av(\Top(T)),
\]
where subtraction is taken in~$\R$.
Thus $R\mapsto R+\Av(R)$ is injective on $\Rp - \{0\}$.
\end{lemma}
\begin{proof}
If $\height(R)=0$, then $\Av(R)=0$. Otherwise its height is strictly
smaller than that of~$R$. In both cases $\Top(T)=\Top(R)$, and
\eqref{E:Av-scale} gives $\Av(\Top(T))=\Av(R)$.
\end{proof}

\begin{proof}[First proof of Theorem~\ref{T:main}]
By Lemma~\ref{L:start}, we may assume that $A$ and $B$ have positive
height. Put $R=\Av(A)$ and $S=\Av(B)$. These are nonzero by
Lemma~\ref{L:Av-height}, and Lemma~\ref{L:Av-rules} gives
\[
 AR=\Av(A^A)=\Av(B^B)=BS.
\]
This common element is nonzero, since a~product of nonzero elements
of $\Rp$ has positive coefficients and no cancellation occurs.
Applying $\Av$ to it yields
\[
 R+\Av(R)=\Av(AR)=\Av(BS)=S+\Av(S).
\]
Lemma~\ref{L:Av-inject} gives $R=S$. Cancel this nonzero factor
by Lemma~\ref{L:cancel} to obtain equality of the representations
of $A$ and $B$, hence $A\cong B$.
\end{proof}

\Needspace{12\baselineskip}
\section{Second proof: three invariants}\label{S:invariants}

This is the proof obtained from ChatGPT, as described in
Section~\ref{S:intro}.

\subsection{The interval sum}
Let $\M(X)$ be the disjoint sum of a~separate copy of every interval
in $\MaxInt(X)$. Even overlapping or isomorphic intervals contribute
separate summands. Each component is bounded and has height
$\height(X)$; in particular, $\comp(\M(X))=|\MaxInt(X)|$.
From the interval correspondence for products and the effect of duality,
\begin{equation}\label{E:M}
\begin{gathered}
 \M(XY)\cong\M(X)\M(Y),\qquad
 \M(X\dual)\cong\M(X)\dual,\\
 \M(nX)\cong n\M(X)\quad(n\geq1).
\end{gathered}
\end{equation}
If $\height(Y)<\height(X)$, then $\M(X+Y)\cong\M(X)$.
The same conclusion holds if the second summand is empty.
All these constructions preserve isomorphisms.

\begin{lemma}\label{L:M-power}
For ordered sets $X$ and $Y$, $\M(X^Y)\cong\M(X)^Y$.
\end{lemma}
\begin{proof}
Write $Y=Y_1+\cdots+Y_r$, with each $Y_i$ connected.
Lemma~\ref{L:intervals} gives
\[
 \M(X^Y)\cong\sum_{(I_i)}\prod_iI_i^{Y_i}
          \cong\prod_i\sum_I I^{Y_i},
\]
where each $I$ ranges over $\MaxInt(X)$.
A~map from $Y_i$ into $\M(X)$ has its image in one summand,
and maps into different summands are incomparable pointwise.
Hence $\sum_I I^{Y_i}\cong\M(X)^{Y_i}$. Finally,
$\prod_i\M(X)^{Y_i}\cong\M(X)^Y$.
\end{proof}

\subsection{Disconnected bases}
\begin{lemma}\label{L:J-power}
Let $P$ be componentwise bounded with $m$ components, and let $Y$
have $c$ components. Then
\[
 \J(P^Y)\cong m^{c-1}\bigl(\J(P)\times Y\dual\bigr).
\]
\end{lemma}
\begin{proof}
Write $P=\sum_{i=1}^mB_i$ and $Y=\sum_{k=1}^cY_k$ as sums of
connected components. The maps $Y\to P$ split according to their
target component on each $Y_k$. Distinct choices give incomparable
suborders, so
\[
 P^Y\cong\sum_{\sigma:\{1,\ldots,c\}\to\{1,\ldots,m\}}
                 \prod_{k=1}^cB_{\sigma(k)}^{Y_k}.
\]
Each factor in a~product is bounded and connected. By
Lemmas~\ref{L:log} and~\ref{L:J-product}, the summand indexed by
$\sigma$ contributes $\sum_k\J(B_{\sigma(k)})\times Y_k\dual$ to
$\J(P^Y)$. For fixed $k$ and $i$, exactly $m^{c-1}$ choices satisfy
$\sigma(k)=i$. Summing the contributions proves the formula.
\end{proof}

\subsection{The invariants}
Let $A$ have height $h\geq1$ and $c$ components, and put
\begin{align*}
P&= \M(A),\\ 
m&= \comp(P),\\ 
\mu&= m^{c-1},\\
W&= \M(\J(P)),\\ 
k&= \comp(W),\\ 
g&= \height(\J(P)).
\end{align*}
Every component of $P$ is bounded, nontrivial, and of height~$h$.
Lemma~\ref{L:J-height} shows that $\J(P)$ is nonempty and
$0\leq g<h$. Thus $W$ is defined and nonempty, its components are
bounded of height~$g$, and $m,k,\mu$ are positive integers.

\begin{proposition}\label{P:invariants}
For $E=A^A$, the constructions
\begin{align*}
 U&= \J(\M(E)),\\ 
 V&= \M(U),\\ 
 Q&= \M(\Js(V))
\end{align*} 
are defined and nonempty. They satisfy
\begin{align}
 U&\cong\mu\bigl(\J(P)\times A\dual\bigr),\label{E:U}\\
 V&\cong\mu\bigl(W\times P\dual\bigr),\label{E:V}\\
 Q&\cong\mu k\,W\dual,\label{E:Q}\\
 h&=\height(U)-\height(Q).\label{E:h-recover}
\end{align}
\end{proposition}
\begin{proof}
Lemmas~\ref{L:M-power} and~\ref{L:J-power} give~\eqref{E:U},
so $U$ is nonempty. Equation~\eqref{E:M} then gives~\eqref{E:V}.
The orders $W$ and $P\dual$ are componentwise bounded.
Apply the upper-cover version of~\eqref{E:J-prod} and use
$\Js(P\dual)=\J(P)\dual$ to obtain
\[
 \Js(V)\cong\mu\bigl(m\Js(W)+k\J(P)\dual\bigr).
\]
The second summand is nonempty and has height~$g$.
If $g>0$, Lemma~\ref{L:J-height}, applied in each component of $W$,
shows that the first summand has height less than~$g$.
If $g=0$, every component of $W$ is a~singleton and $\Js(W)$ is empty.
In either case all intervals of maximum height come from the second
summand. Thus $Q$ is defined and
\[
 Q\cong\mu k\M(\J(P)\dual)\cong\mu kW\dual,
\]
which is nonempty. Finally, $\height(U)=g+h$ and $\height(Q)=g$,
proving~\eqref{E:h-recover}.
\end{proof}

\begin{proof}[Second proof of Theorem~\ref{T:main}]
By Lemma~\ref{L:start}, assume that $A$ and $B$ have positive height.
Use subscripts $A$ and $B$ for the quantities just defined.
The constructions $U,V,Q$ depend only on $E\cong A^A\cong B^B$,
so the two representations give isomorphic results.
Fix a~connected ordered set~$Z$, and write
$[X]=\Hom(Z,X)$ and $[X]'=\Hom(Z\dual,X)$.

Apply $\Hom(Z\dual,-)$ to~\eqref{E:Q} and $\Hom(Z,-)$
to~\eqref{E:V}. Equations~\eqref{E:hom2} and~\eqref{E:hom3} give
\begin{align*}
 \mu_Ak_A[W_A]=\mu_Bk_B[W_B],\\
 \mu_A[W_A][P_A]'=\mu_B[W_B][P_B]'.
\end{align*}
The counts and scalar factors are positive. Dividing the second
equality by the first gives
\begin{align*}
 \frac{[P_A]'}{k_A}&= \frac{[P_B]'}{k_B},\\
  k_B[P_A]'&= k_A[P_B]'.
\end{align*}
As $Z\dual$ ranges over the connected ordered sets,
Lemma~\ref{L:lovasz} gives $k_BP_A\cong k_AP_B$.
Applying $\J$, which commutes with disjoint sums, gives
\[
 k_B\J(P_A)\cong k_A\J(P_B).
\]
Put $a=[\J(P_A)]>0$ and $b=[\J(P_B)]>0$.
The last isomorphism gives $k_Ba=k_Ab$, while~\eqref{E:U} gives
\[
 \mu_Aa[A\dual]=\mu_Bb[B\dual].
\]
Substitute $b=k_Ba/k_A$ and cancel $a$ to obtain
\[
 \mu_Ak_A[A\dual]=\mu_Bk_B[B\dual].
\]
Taking $Z=\one$ gives $\mu_Ak_A|A|=\mu_Bk_B|B|$.
Equation~\eqref{E:h-recover} independently gives $h_A=h_B>0$,
so Lemma~\ref{L:start} yields $|A|=|B|$.
It follows that $\mu_Ak_A=\mu_Bk_B$, and hence
$[A\dual]=[B\dual]$ for every connected~$Z$.
A~final application of Lemma~\ref{L:lovasz} gives
$A\dual\cong B\dual$, and therefore $A\cong B$.
\end{proof}

\section{Comparison with the earlier proof}\label{S:compare}

The proof in~\cite{FG} uses Hashimoto's refinement
theorem~\cite{HasEA}, the canonical form of
McKenzie~\cite{McKJJ}, and Farley~\cite{FarBC} to define a~logarithm
on connected ordered sets. It also constructs an ordering of connected
isomorphism types compatible with products and connected exponentiation.
That ordering selects one interval of maximum height, from which
the reconstruction proceeds. The same constructions treat a~family
of iterated exponential expressions.

The present proofs require the logarithm only on bounded orders.
The construction $\Log(B)=\J(B)\dual$ supplies it without direct
factorization or canonical coordinates. Working with every interval
of maximum height removes the need to choose an auxiliary ordering.
Homomorphism counts supply cancellation, and Lemma~\ref{L:lovasz}
includes a~proof of the counting result itself.

\Needspace{17\baselineskip}
The table records the ingredients used in the three arguments.

\begin{center}\small
\begin{tabular}{@{}>{\raggedright\arraybackslash}p{0.44\textwidth}ccc@{}}
\toprule
Ingredient & Earlier~\cite{FG} & Average & Invariants\\
\midrule
Intervals of maximum height & $\checkmark$ & $\checkmark$ & $\checkmark$\\
Hashimoto's refinement theorem & $\checkmark$ & -- & --\\
McKenzie--Farley canonical form & $\checkmark$ & -- & --\\
Auxiliary ordering of isomorphism types & $\checkmark$ & -- & --\\
Elements with exactly one lower cover & -- & $\checkmark$ & $\checkmark$\\
Homomorphism counts & -- & $\checkmark$ & $\checkmark$\\
Iterated exponentiation & $\checkmark$ & -- & --\\
\bottomrule
\end{tabular}
\end{center}

The proofs share the interval lemma and the lower-cover construction.
The first short proof packages the intervals into a~rational average and uses the strict height decrease to recover its unknown factor. 
The second  short proof
retains integral multiplicities and uses the three invariants $U,V,Q$
to compare homomorphism counts. Neither argument requires the
reconstruction theorem or the auxiliary ordering from~\cite{FG}.

\Needspace{15\baselineskip}

\end{document}